\documentclass[10pt,a4paper,twoside]{amsart}
\usepackage[T1]{fontenc}
\usepackage[utf8]{inputenc}
\usepackage[english]{babel}
\usepackage{amsmath}
\usepackage{amsthm}
\usepackage{amssymb}
\usepackage{amsfonts}
\usepackage{amsxtra}
\usepackage{enumerate}
\usepackage{verbatim}
\usepackage{color}
\usepackage[mathscr]{eucal}

\newcommand{\N}{\mathbb N}
\newcommand{\Z}{\mathbb Z}

\newcommand{\R}{\mathbb R}

\newtheorem{satz}{Satz}[section]

\newtheorem*{theo}{Theorem}
\newtheorem{theorem}[satz]{Theorem}
\newtheorem{definition}[satz]{Definition}
\newtheorem{lemma}[satz]{Lemma}
\newtheorem{cor}[satz]{Corollary}

\newtheorem{prop}[satz]{Proposition}
\newtheorem*{remark}{Remark}
\newtheorem{exmp}[satz]{Example}

\DeclareMathOperator*{\conv}{conv}

\DeclareMathOperator*{\dist}{dist}

\DeclareMathOperator*{\inte}{int}

\DeclareMathOperator*{\Light}{Light}

\let\phi=\varphi
\let\e=\varepsilon

\begin{document}

\setcounter{page}{1}

\title{Closed geodesics in Lorentzian surfaces}

\author[Stefan Suhr]{Stefan Suhr}
\address{Fachbereich Mathematik, Universit\"at Hamburg}
\email{stefan.suhr@mathematik.uni-hamburg.de}

\date{\today}

\begin{abstract}
We show that every closed Lorentzian surface contains at least two closed geodesics. Explicit examples show the optimality of this claim. Refining this result we 
relate the least number of closed geodesics to the causal structure of the surface and the homotopy type of the Lorentzian metric. 
\end{abstract}
\maketitle

\section{Introduction}\label{last}

The\footnote{53C22, 53C50, closed geodesics, Lorentzian manifolds} study of closed geodesics in Riemannian geometry has been the source of deep insights 
(e.g. Morse theory) into the properties of positive definite variational problems and the structure of Riemannain manifolds in general. A natural attempt is thus to 
transfer this body of knowledge to Lorentzian and pseudo-Riemannian manifolds. Some of the topological techniques known from Riemannian geometry carry over 
only in rather special cases, namely in the presence of a (timelike) Killing vector field (compare \cite{mas1} or \cite{bmp1}). In general all Morse theoretic tools fail, 
though. This is due to the fact that the index of a geodesic in a pseudo-Riemannian manifold is always infinite. The known results addressing more generic 
situations, i.e. without Killing vector fields, rely on maximizing techniques for the length functional. The classical theorems of Tipler (\cite{tip1}) and Galloway 
(\cite{ga1}) fall into this category. 

The best result, with respect to generality, is the theorem by Galloway (\cite{ga2}) stating that any closed, i.e. compact without boundary, Lorentzian surface 
contains at least one closed timelike or lightlike periodic geodesic. Here we propose the following distinction between closed and periodic lightlike geodesics. A 
lightlike geodesics $\gamma \colon I\to M$ in a pseudo-Riemannian manifold $(M,g)$ is periodic if there exists $s,t\in I$ with $\gamma(s)=\gamma(t)$ and 
$\dot{\gamma}(t)$ is collinear to  $\dot{\gamma}(s)$. On the other hand call $\gamma$ closed if $\dot{\gamma}(t)=\dot{\gamma}(s)$, i.e. the orbit of the geodesic 
flow of $(M,g)$ defined by $\dot{\gamma}$ is closed. The difference between closed and periodic but non-closed can be characterized in terms of completeness: 
A lightlike periodic geodesic $\gamma$ is closed iff $\gamma$ is complete (compare \cite{ca1} for a discussion). The proof relies on the observation that each 
cycles of a periodic geodesic is a reparameterization of the previous cycle by an affine function. If the slope of this affine function is larger than $1$, the 
reparameterization accelerates the geodesic and the geodsic will be incomplete in the positive direction. If the slope of the affine function lies strictly between $0$ 
and $1$, the geodesic is decelerated anf the geodesic is incomplete in the negative direction. Another noticeable difference between closed and periodic 
but non-closed geodesics is the fact that the later ones are not critical points of the energy functional on the space of loops. This is an easy consequence of the 
formula for the first variation of the energy. Refer to example \ref{E58} and the example \ref{E60} for periodic and  non-closed lightlike geodesics.

We will explain in example \ref{E58} why the argument in the proof of Galloway's theorem implies the existence of a closed timelike or periodic lightlike geodesic
only. Therefore it is not clear at this point if the geodesic flow of a closed Lorentzian surface contains any closed flowlines. These notes will settle this question with 
the following theorem (see section \ref{last3}):
\begin{theo}
Every closed Lorentzian surface contains two closed geodesics, one of which is definite, i.e. time- or spacelike. 
This lower bound is optimal.
\end{theo}

The least number of closed geodesic in a Lorentzian surface $(M,g)$ depends in a precise way on the connected component 
of $g$ in the set of Lorentzian metrics (proposition \ref{T61}). Explicit examples show this relation to be optimal
(example \ref{E60}). 

The text is structured as follows. First we will review some results concerning the
existence of closed causal geodesics in closed Lorentzian manifolds. The focus will lie mainly on results which restrict  
the causal structure of the Lorentzian manifold only. This will leave out the entire literature on closed geodesics in 
stationary or static space-times. 

In section \ref{last3} we study the problem of the least number of closed geodesics in a closed pseudo-Riemannian (i.e. Lorentzian) surface. We will show that 
every closed Lorentzian surface contains a timelike or spacelike closed geodesic. Further the geodesic flow of a closed Lorentzian surface contains at least two 
closed flowlines. Previously known examples show that these lower bounds on the multiplicity of closed geodesics are optimal. For a spacetime structure on the 
torus $T^2$ we are able to improve this number to four. The methods to derive these results rely on an analysis of the space of closed curves with constant causal 
character.

\subsection*{Notation}

$M$ always denotes a closed surface of vanishing Euler charateristic, i.e. $M$ is diffeomorphic to the $2$-torus $T^2$ or the Klein bottle $K^2$.
We always assume that a fixed complete Riemannian metric $g_R$ is chosen on all manifolds. If not noted otherwise Riemannian 
metrics on covering spaces are assumed to be equal to the lifted metric.

By a {\it nonspacelike} curve we understand the usual notion of piecewise $C^1$ nonspacelike curves (see \cite{onei} pg. 146). A curve is {\it nontimelike} if it is 
nonspacelike for $(M,-g)$. A curve is {\it causally constant} if it is either nonspacelike or nontimelike. 

By $\Lambda M_\pi$ we understand the space of all continuous loops representing the free homotopy class $\pi$, 
and by $\Lambda(M,[g])_\pi$ we understand the space of causally constant loops representing $\pi$, where $[g]$ denotes 
the conformal class of $g$.

\section{Historical remarks on closed causal geodesics}\label{HR}

The study of closed causal geodesics in Lorentzian geometry was initiated by F. Tipler in \cite{tip1} with the following result.
\begin{theorem}[\cite{tip1}]\label{TTIP}
Let $(M,g)$ be a compact spacetime with a covering space containing a compact Cauchy hypersurface. Then $(M,g)$ 
contains a closed timelike geodesic.
\end{theorem}
The method used in the proof is an adaptation of the Birkhoff process, known in Riemannian geometry, to minimize
the arclength in $\Lambda M_\pi$. It relies heavily on the 
causal structure of the covering space. The last development in this direction is due to Guediri (\cite{guediri1}) 
and, with a generalized version, Sanchez in \cite{sa2}.
\begin{prop}[\cite{sa2}]\label{P50}
Let $(M,g)$ be a compact Lorentzian manifold admitting a globally hyperbolic cover $(\overline{M},\overline{g})$. 
Assume that a conjugacy class $\mathcal{C}\subset \mathcal{D}(\overline{M},M)$ satisfies:
\begin{enumerate}
\item $\mathcal{C}$ contains a closed timelike curve $\gamma$.
\item $\mathcal{C}$ is finite.
\end{enumerate}
Then there exists at least one closed timelike geodesic in $\mathcal{C}$.
\end{prop}

Galloway formalized the adaptation of the Birkhoff process by introducing $t$-homotopies.
\begin{definition}[\cite{ga1}]
Let $(M,g)$ be a Lorentzian manifold. 
A continuous mapping $H\colon I \times [0,1]\to M$, subject to certain boundary conditions (either fixed endpoints or 
closed curves) is called a t-homotopy if the curves $t\mapsto H(t,\tau)$ are causal for all $\tau \in [0,1]$. 
\end{definition}
The notion of $t$-homotopy reformulates naturally to a statement about the connected components of the set of causal 
curves in the set of 
continuous mappings $I\to M$.  The Riemannian metric $g_R$ introduces a complete metric topology on $C^0(I,M)$. 
Two causal curves $\gamma_{1,2}\colon I\to M$ are then $t$-homotopic
iff they lie in the same connected component of the set of causal curves with domain $I$ (this is always understood with 
respect to given boundary conditions and therefore readily extend to, for example, loops). $t$-homotopy classes
are therefore the connected components of the set of causal curves. 

To state the condition for the existence of closed timelike geodesics, recall the definition of wideness from \cite{ga1}: Given two Lorentzian metrics $g_1,g_2$ on 
$M$, $g_2$ is said to be {\it wider} than $g_1$ if every $g_1$-causal vector is $g_2$-timelike. This algebraic condition translates to the geometric picture of the 
time cones of $g_2$ being larger than those of $g_1$. Drawing inspiration from globally hyperbolic spacetimes, Galloway poses the following condition: A 
$t$-homotopy class $\mathcal{C}\subset\Lambda M$ of a given Lorentzian metric $g$ is {\it stable} if there exists a Lorentzian metric $\hat{g}$, wider than $g$, 
such that 
$$\sup_{\gamma\in \mathcal{C}}L^{\hat{g}}(\gamma)<\infty,$$
where $L^{\hat{g}}$ denotes the length functional of $\hat{g}$. Galloway then goes on to prove the following theorem.

\begin{theorem}[\cite{ga1}]\label{T39}
Let $(M,g)$ be a compact Lorentzian manifold. Then each stable free $t$-homotopy class $\mathcal{C}\subset\Lambda M$ 
contains an arclength maximizing closed causal curve. If $\sup_{\gamma\in\mathcal{C}}L^g(\gamma)>0$, this curve 
is necessarily a closed timelike geodesic.
\end{theorem}

The proof consists of a reformulation of the Birkhoff's shortening process for loops in Riemannian manifolds
for the case of Lorentzian manifolds. 

In a subsequent study Galloway then proved the following result (Note that we use the distinction introduced in section \ref{last}).
\begin{theorem}[\cite{ga2}]\label{T40}
Every closed Lorentzian surface contains a closed timelike  or periodic lightlike geodesic.
\end{theorem}

The proof of theorem \ref{T40}, given in \cite{ga2}, does not imply the existence of a closed geodesic, i.e. the tangent curve is a periodic orbit of the geodesic flow, 
in any closed Lorentzian surface. To make this claim precise we will show that the construction in \cite{ga2} does not necessarily yield a complete periodic (i.e. 
closed) lightlike geodesic. 

First we give a sketch of Galloway's argument and then discuss an example displaying the claimed phenomenon.

\begin{proof}[Sketch of proof of theorem \ref{T40}]
Let $M$ be a closed $2$-manifold equipped with a Loren\-tzian metric $g$. Since periodic geodesics in any covering manifold $M'\to M$ project to periodic 
geodesics in $M$, we can assume that $M$ is orientable and time orientable. This together with the assumption that dim $M=2$ implies that $M$ is 
diffeomorphic to a $2$-torus.

One of the special features of Lorentzian metrics $g$ on $2$-manifolds is, that $-g$ is Lorentzian as well. Since it is well known that compact Lorentzian manifolds 
contain closed timelike curves, it follows that $(M,g)$ contains a smooth closed spacelike curve, i.e. an immersed compact smooth spacelike hypersurface 
$\Sigma$. We can even choose $\Sigma$ to be embedded. i.e. simple closed. 

It follows now that $\Sigma$ does not separate $M$, since this would yield a contractible spacelike curve in the universal cover.
Therefore it is possible that $\Sigma$ is not acausal (in fact this will be the case in our subsequent example).

This problem can be resolved by passing to the so-called Geroch-covering $M_\Sigma$ of $M$ relative to $\Sigma$ (\cite{ger0}). The Geroch covering of a 
spacetime $(M,g)$ relative to a complete spacelike hypersurface $\Sigma$ is defined as follows: The time orientation of $(M,g)$ induces an orientation of the 
normal bundle $T^\perp \Sigma$ of $\Sigma$. Fix a point $p_0\in M$ and consider the set 
$$C:=\{(p,\gamma)|\;p\in M, \gamma \text{ is a continuous path from $p$ to $p_0$}\}.$$
Define the equivalene relation $``\sim"$ by $(p,\gamma)\sim (p',\gamma')$ if $p=p'$ and the intersection numbers of $\gamma$ and $\gamma'$ with $\Sigma$ 
coincide. The manifold $M_\Sigma:=C/\sim$ is called the Geroch covering of $M$ relative to $\Sigma$. 

$\pi_\Sigma\colon M_\Sigma\to M$, $[(p,\gamma)]\mapsto p$ is a regular cover such that (1) each component $\widetilde{\Sigma}$ of $\pi_\Sigma^{-1}(\Sigma)$ is 
a compact sapcelike hypersurface 
which separates $M_\Sigma$ and (2) $\overline{J^-(p)}\cap J^+(\widetilde{\Sigma})$ as well as $\overline{J^+(p)}\cap J^-(\widetilde{\Sigma})$ are compact for all 
$p\in M_\Sigma$. 

By property (1) every component $\widetilde{\Sigma}$ is acausal. Now if $\widetilde{\Sigma}$ is a Cauchy-hypersurface the theorem of Tipler (theorem \ref{TTIP})
implies the existence of a closed timelike geodesic and we are done. 

If $\widetilde{\Sigma}$ is not a Cauchy-hypersurface the Cauchy-horizon $H(\widetilde{\Sigma})$ of $\widetilde{\Sigma}$ is not empty (see \cite{ger1} for the 
definition). Since the argument is invariant under switching the time orientation, we can assume that $H^+(\widetilde{\Sigma})\neq\emptyset$. The remainder of the 
proof now consists of showing that any past inextendible null geodesic generator of $H^+(\widetilde{\Sigma})$ projects to a periodic lightlike geodesic in $M$. 
\end{proof}

Now we will show via an example that the geodesics constructed above are not necessarily closed. Notice that these geodesics are nevertheless periodic in
our terminology. The following family of metrics is a slight modification of the first example given in \cite{ga2} (see below).
\begin{exmp}\label{E58}
Denote with $x,y$ the standard coordinates on $\R^2$. Consider on $\R^2$ the Lorentzian metric 
$$\overline{g}_\e:= (\cos^2 (2\pi x) -\e) [dy^2-dx^2] -2\sin(2\pi x) dxdy$$
where $\e\in (0,1/4)$. Since $\R^2$ is simply connected, $(\R^2,\overline{g}_\e)$ is time orientable. Choose the time orientation such that $\partial_x$ is 
futurepointing in $x=0$.  Observe that $(\R^2,\overline{g}_\e)$ is not globally hyperbolic. This follows directly from the fact that $J^+((0,0))\cap J^-((0,2\pi))$ contains
$\{\frac{\pi}{2}\}\times (-\infty,0)$ and is therefore not compact.

It is obvious that $\overline{g}_\e$ descends to a Lorentzian metric $g_\e$ on the $2$-torus $\R^2/\Z^2$. Observe that $g_\e$ is time orientable. 
The set $\{x=0\}$ projects to a smooth compact spacelike hypersurface $\Sigma\subseteq T^2:=\R^2/\Z^2$. $\Sigma$ is not acausal in $(T^2,g_\e)$  by 
construction. Observe that in every point $p\in \R^2$ there exists a $v_y\in \R$ such that $\partial_x+v_y\partial_y$ is timelike future pointing. Since 
$\overline{g}_\e$ is $\Z^2$-invariant this implies that there exists a timelike curve $\gamma\colon [0,1]\to \R^2$ with $\gamma(0)\in \{x=0\}$ and $\gamma(1)\in 
\{x=1\}$. Projecting $\gamma$ to $T^2$ yields a chronological curve with endpoints on $\Sigma$. 

According to \cite{sa1} (or by direct calculation of the geodesic equations) the lightlike geodesics parameterizing $\{\cos^2 2\pi x =\e\}$ are not complete. By 
construction these geodesics project to periodic lightlike geodesics in $T^2$. 

It is easy to see that the Geroch covering $T^2_\Sigma$ of $T^2$ relative to $\Sigma$ is $\R^2/[\{0\}\times \Z]$. 
Notice that the induced metric is not globally hyperbolic, since the covering space $(\R^2,\overline{g}_\e)$ is not globally hyperbolic. It is then easy to see that the 
Cauchy horizon of $\widetilde{\Sigma}:=\{x=0\}\subset T^2_\Sigma$ is equal to $\{x=\pm x_\e\}$ where $x_\e$ is the smallest positive real number such that 
$\cos^2 2\pi x_\e =\e$. We have seen above that the null geodesic generators of $H(\widetilde{\Sigma})$ are incomplete. Therefore they are periodic, but not 
closed lightlike geodesics. In fact this is the case for all periodic lightlike geodesics in $(T^2,g_\e)$.

One should note, though, that $(T^2,g_\e)$ contains closed timelike geodesics (e.g. $\{x=1/4\}$), but these do not intersect $D(\Sigma)\cup H(\Sigma)$. 
\end{exmp}

The original example in \cite{ga2} is the metric $-\overline{g}_\e$ with $\e=0$. For later purposes we explain it in detail.
\begin{exmp}\label{E59}
Consider $\R^2$ with the natural coordinates $\{x,y\}$ and the Lorentzian metric
$$\overline{g}:=\cos^2(2\pi x)[dx^2-dy^2]+2\sin(2\pi x)dxdy.$$ 
$\overline{g}$ is obviously invariant under the translations $(x,y)\mapsto (x,y+t)$ for all $t\in\R$ and the map $(x,y)\mapsto (x+1/2,-y)$. Thus the metric 
$\overline{g}$ descends to a Lorentzian metric $g$ on the Klein bottle $K^2=\R^2/\Gamma$ where $\Gamma:=\langle (x,y)\mapsto (x,y+1), (x,y)\mapsto 
(x+1/2, -y)\rangle$. In \cite{ga2} it is shown that $(K^2,g)$ does not contain any spacelike closed geodesics. 

Using \cite{sa1} we see that the geodesic parameterization of the unique closed smooth lightlike curve in $(K^2,g)$ is complete. Therefore it is a closed lightlike 
geodesic in our terminology. Again from \cite{sa1} follows that $[x=0]$ is the trace of the single closed timelike geodesic of $(K^2,g)$. Concluding we see that 
$(K^2,g)$ contains exactly two geometrically distinct closed geodesics.
\end{exmp}

Finally we mention two results concerning closed timelike geodesics in Lorentzian manifolds satisfying special assumptions.
\begin{prop}[\cite{guediri1}]
Let $(T^2,g)$ be a Lorentzian $2$-torus with geodesically connected universal covering. Then, it contains a closed 
timelike geodesic.
\end{prop}

\begin{theorem}[\cite{FJP}]
Let $(M,g)$ be a compact Lorentzian manifold with dim $M\ge 2$ that admits a Killing vector field $K$ that is timelike somewhere. Then there is some non trivial 
periodic non self-intersecting timelike geodesic in $(M, g)$. If either one of the following two conditions is satisfied, then there are at least two non trivial periodic 
non self-intersecting geodesics in $M$:
\begin{itemize}
\item[(a)] $\max_{q\in M} g(K_q,K_q)\neq 0$;
\item[(b)] $K$ is never vanishing.
\end{itemize}
When either condition is satisfied, if in addition $K$ has at most one periodic integral line, then there are infinitely many geometrically distinct non trivial periodic 
non self-intersecting geodesics in $(M, g)$.
\end{theorem}

It should be pointed out, that example \ref{E59} shows the optimality of one part of the theorem. The part claiming that any Lorentzian manifold with a Killing vector field, 
timelike somewhere and subject to the conditions (a) or (b), contains at least two non trivial periodic non self-intersecting closed geodesics is optimal by the example.

\section{Generalities about Lorentzian metrics}

\subsection{The lightlike distributions}
Locally every Lorentzian surface $(M,g)$ gives rise to two transversal lightlike distributions. In general these distributions are not globally well defined. Note that 
they are globally well defined if and only if the underlying manifold is orientable. This can be seen as follows: If the lightlike distributions are well defined, every 
choice of $g_R$-unit vector field lying in one of the distributions is well defined up to multiplication with $-1$. Note that the sign of 
$g|_{\conv(v,w)\times\conv(v,w)}$ is locally constant for any continuous choice of lightlike basis $\{v,w\}$ of the tangent spaces. In the case of a $2$-dimensional 
vector space the bases $\{v,w\}$ and $\{-v,-w\}$ share the same orientation. Therefore $M$ is orientable. Conversely if $M$ is orientable, we can locally choose a 
positively oriented basis in the lightlike distributions. Again this basis is globally well defined only up to sign. But in this case the distributions are still well defined.

Assume now that $M$ is orientable. Then $M\cong T^2$ and $\pi_1(M)\cong H_1(M,\Z)\cong \Z^2$ operates by a co-compact 
action on the universal cover $\widetilde{M}$. Furthermore the image of any pair of generators $v,w$ of $H_1(M,\Z)$ 
defines a basis of $H_1(M,\R)$. Choose such a base $\{v,w\}\subset H_1(M,\Z)$ of $H_1(M,\R)$ and let 
$\{\alpha_1,\alpha_2\}$ be the dual basis of 
$H^1(M,\R)$ with representatives $\omega_i\in\alpha_i$ (We will always use the deRham-cohomology, i.e. $\omega_i$ is a 
smooth closed $1$-form). For an absolutely continuous curve $\zeta\colon [a,b]\to M$ define $\zeta(b)-\zeta(a)\in H_1(M,\R)$ 
by the equation
$$\alpha_i(\zeta(b)-\zeta(a)):=\int_\zeta \omega_i.$$
This definition depends of course on the chosen representatives $\omega_i\in\alpha_i$, but for any pair of representatives there 
exists a constant such that the difference between both of $\zeta(b)-\zeta(a)$ is uniformly bounded independent of $\zeta$.

We associate to a lightlike distribution $\mathfrak{D}$ on $M$ the class $m^{\mathfrak{D}}\in P H_1(M,\R)$ in the projective 
space over the first real homology vector space of $M$ as follows: Consider a curve $\zeta\colon \R \to M$ everywhere 
tangential to $\mathfrak{D}$. Then there exists a unique line $l\subset H_1(M,\R)$ such that the distance of any 
$\zeta(T_2)-\zeta(T_1)$ to $l$ is bounded by a uniform constant depending only on $\mathfrak{D}$. 
The line $l$ is independent of $\zeta$. This follows from \cite{hehi} Part A, section 4.3, since
$\zeta$ parameterizes (not necessarily injective) a leaf of a nonsingular foliation of $M$. 
Call $m^\mathfrak{D}:=[l]\in PH_1(M,\R)$ the {\it rotation class} of $\mathfrak{D}$. Introduce an arbitrary norm $\|.\|$ on 
$H_1(M,\R)$. Note that the rotation class satisfies
$$m^\mathfrak{D}=\lim_{\|\zeta(T)-\zeta(T')\|\to\infty}[span(\zeta(T)-\zeta(T'))]\in PH_1(M,\R)$$
for any curve $\zeta$ tangential to $\mathfrak{D}$.

Call $(M,g)$ space orientable if $(M,g)$ admits a spacelike nonsingular vector field. Note that this is equivalent to 
$(M,-g)$ being time orientable. In this notation the following conditions are equivalent:
\begin{enumerate}[(i)]
\item The lightlike distributions are orientable.
\item $(M,g)$ is time and space orientable.
\item $M$ is orientable and $(M,g)$ is time orientable
\item $M$ is orientable and $(M,g)$ is space orientable
\end{enumerate}
Recall that any Lorentzian manifold admits a twofold time orientable covering. Therefore any 
Lorentzian manifold admits a, at most, fourfold orientable and time orientable covering. 

Assume now that the lightlike distributions are well defined and orientable, i.e. there exists two future pointing lightlike 
vector fields $X^+$ and $X^-$ such that $\{X^+_p,X^-_p\}$ is a positively oriented basis of $TM_p$ for all $p\in M$.
Define $\mathfrak{D}^+$ through $X^+\in \mathfrak{D}^+$ and $\mathfrak{D}^-$ through $X^-\in \mathfrak{D}^-$. 
Abridge $m^\pm:=m^{\mathfrak{D}^\pm}$.

\subsection{Lorentzian metrics on surfaces}

We collect in this paragraph some more or less well known facts about the connected 
components of the set of Lorentzian metrics in a suitable geometric way. 

A vector field without zeros is said to be causally constant, if either $g(X,X)\ge 0$ or $g(X,X)\le 0$ everywhere. 
The existence of causally constant vector fields is equivalent to space resp. time orientability. 

Assume $(M,g)$ to be space and time orientable. Then $M$ is diffeomorphic to the $2$-torus $T^2=\R^2/\Z^2$ and 
$TM$ is parallelizable.  Consider the space of unit vector fields $\Gamma(T^{1,R}M)$ of $(M,g_R)$.
Consider then the projection $\pi_{S^1}\colon  T^{1,R}M\to S^1$ defined by the 
composition of a bundle isomorphism $T^{1,R}M\cong M\times S^1$, where the isomorphism is induced by the differential 
of any diffeomorphism $S^1\times S^1\to M$, with the projection onto the second factor of 
$ M\times S^1$. Thus for a unit vector field $X\in \Gamma(T^{1,R}M)$ and a closed curve $\gamma\colon S^1\to M$, 
the degree of $X\circ \gamma$ is naturally defined as the degree of $\pi_{S^1}\circ X\circ \gamma$. Fixing a closed curve 
$\gamma\colon S^1\to M$ the degree of $X\circ \gamma$ is independent of the choice of causally constant vector 
field $X\in \Gamma(T^{1,R}M)$.

\begin{definition}
Let $(M,g)$ be space and time orientable. For $\sigma\in \pi_1(M)$ define the {\it rotation number} 
$n_\sigma(g)$ of $g$ along $\sigma$ to be the degree of $X\circ\gamma$, where $X$ is any causally constant $g_R$-unit 
vector field and $\gamma$ is any curve representing $\sigma$.
\end{definition}

Note that for a nonspacelike or nontimelike curve $\gamma$ the rotation number of $g$ along the homotopy class 
$[\gamma]$ vanishes. Then the map $\Gamma(T^{1,R}M)\to Lor(M)$, $X\mapsto g_R-2X^\flat\otimes X^\flat$ 
$(\flat:=\flat^{g_R})$ induces an isomorphism of $\pi_0(\Gamma(T^{1,R}M))$ to $\pi_0(Lor(M))$. Therefore 
two Lorentzian metrics $g$ and $g'$ on $M$ are homotopic in $Lor(M)$ if and only if for a pair of generators $\{\sigma,\tau\}$ 
of $\pi_1(M)$, $n_{\sigma}(g)=n_{\sigma}(g')$ and $n_{\tau}(g)=n_{\tau}(g')$.

\begin{prop}\label{P49}
Let $(M,g),(N,g')$ be space and time orientable closed Lorentzian surfaces.
\begin{enumerate}
\item Set $k_g:=\min\{|n_\sigma(g)|+|n_\tau(g)|\}$, where the 
minimum is taken over all pairs of generators $\{\sigma,\tau\}$ of $\pi_1(M)$. Then $g$ is isometric to a Lorentzian 
metric $g'$ with $n_{\sigma}(g')=k_g$ and $n_{\tau}(g')=0$.
\item  Assume for $(M,g)$ and $(N,g')$, $k_g=k_{g'}$. 
Then there exists a diffeomorphism $F\colon M \to N$ such that $F^\ast g'$ is homotopic to $g$ in $Lor(M)$.
\end{enumerate}
\end{prop}

\begin{proof}
(1) Let $\{\sigma,\tau\}$ be a pair of generators of $\pi_1(M)$ and set $m:=n_\sigma(g)$, $n:=n_\tau(g)$. 
Define $\widetilde{m}:=\frac{m}{gcd(m,n)}$ and $\widetilde{n}:=\frac{n}{gcd(m,n)}$. Choose $a,b\in\mathbb Z$ 
such that $\widetilde{n}b+\widetilde{m}a=1$. If $n_{a\sigma+b\tau}(g)\ge 0$ set $A:=\left(\begin{smallmatrix} 
                       \widetilde{n} & a \\ -\widetilde{m} & b
                    \end{smallmatrix}
                   \right)$. Otherwise define $A:=\left(\begin{smallmatrix} 
                       \widetilde{n} & -a \\ -\widetilde{m} & -b
                    \end{smallmatrix}\right)$.                   
Then there exists a diffeomorphism of $F\colon M\to M$ with $F_\ast \colon \pi_1(M)\to \pi_1(M)$ identical 
to $A$ relative to the base $\{\sigma,\tau\}$. For $g':=F^\ast g$ holds $n_{\sigma}(g')=k_g$ and $n_{\tau}(g')=0$, 
since $n_{\sigma}(g')$ is a generator of the image of $n_{.}(g')\colon \pi_1(M)\to \Z$, $\zeta\mapsto n_\zeta(g)$. 

(2) Taking the previous remarks into account, (2) follows readily from (1). 
\end{proof}

\section{Closed geodesics in Lorentzian surfaces}\label{last3}

Theorem \ref{T40} is unsatisfactory when adopting the point of view of dynamical systems, especially closed orbits of geodsic flows, towards Lorentzian geometry. 
As we have seen above, periodic lightlike geodesics are i.g. not closed geodesics in the sense that the tangent curves are not closed orbits of the geodesic flow. 
This leaves us with the possibility for Lorentzian surfaces without any closed geodesics in the terminology defined in the introduction. To clarify this gap we claim 
the following theorems (We call a geodesic 
$\gamma$ {\it definite} if $g(\dot{\gamma},\dot{\gamma})\neq 0$).
\begin{theorem}\label{T60}
Every closed Lorentzian surface contains a definite (i.e. timelike or spacelike) closed geodesic. This lower bound is optimal.
\end{theorem}

\begin{theorem}\label{T60a}
The geodesic flow of a closed Lorentzian surface has at least two closed orbits. This lower bound is optimal.
\end{theorem}
Example \ref{E59} shows the optimality of the assertions. The induced Lorentzian metric on the Klein bottle has exactly two 
closed geodesics, one is timelike the other lightlike. Note that for Lorentzian surfaces spacelike geodesics locally maximize 
the energy functional among spacelike curves. This is due to the fact that for $dim\, M=2$ and $(M,g)$ 
Lorentzian, $(M,-g)$ is Lorentzian as well. 

The theorems follow from a seperate analysis of two different cases. 

The method of proof for theorem \ref{T60a} shows the following corollary as well.
\begin{cor}
Every time and space orientable closed Lorentzian surface contains four closed geodesics, two of which must be definite, i.e. 
every spacetime structure on the $2$-torus contains at least four closed geodesics.
\end{cor}

\subsection{class A and class B surfaces}

\begin{definition}
$(M,g)$ belongs to class $A$ if there exists a finite time and space orientable cover $(M',g')$ such that $m^+\neq m^-$ 
for $(M',g')$. The complementary case $m^+=m^-$ is denoted with class $B$. 
\end{definition}

Note that this definition is independent of the chosen finite cover.

In the case that $(M,g)$ is class A, time and space orientable, the homology classes $\zeta(T)-\zeta(T')$ lie at a 
bounded distance from two halflines $\overline{m}^+,\overline{m}^-$ of $m^+,m^-$ for $\zeta$ future pointing and 
$T'\le T$. The set $\conv(\overline{m}^+\cup\overline{m}^-)$ is a proper cone in $H_1(M,\R)$. 

We can characterize class A in more traditional terms of causality theory.
\begin{prop}\label{P60}
Let $(M,g)$ be a closed Lorentzian surface. Then the properties are equivalent.
\begin{enumerate}[(i)]
\item $(M,g)$ is class $A$,
\item the time orientation covers of $(M,g)$ and $(M,-g)$ are vicious and 
\item the time orientation cover of $(M,g)$ is vicious and the universal cover is globally hyperbolic.
\end{enumerate}
\end{prop}
Part of the proposition goes back to the diploma thesis of E. Schelling (\cite{schelling}). Satz 4.1 therein 
proves that the Abelian cover of any class A surface is globally hyperbolic (see appendix \ref{AB}).

\begin{proof}
Both viciousness and global hyperbolicity lift to finite coverings and pass to finite factors. Therefore we can 
assume that $(M,g)$ is time and space orientable.

(i)$\Rightarrow$ (ii): Assume that $(M,g)$ is class $A$. It suffices to show viciousness for $(M,g)$, since for $(M,-g)$ 
and the appropiate choice of time orientation of $(M,-g)$, the set $\{X^-_p,-X^+_p\}$ is a positively oriented basis of future
pointing vectors. Therefore the rotation halflines for $(M,-g)$ are $\overline{m}^-$ and $-\overline{m}^+$. 

Denote with $\mathfrak{F}^\pm$ the foliations of $M$ tangential to $\mathfrak{D}^\pm$ and with $\widetilde{\mathfrak{F}}^\pm$ 
their lift to the universal cover $\widetilde{M}$. 
Every leaf of the foliations $\widetilde{\mathfrak{F}}^\pm$  separates $\widetilde{M}\cong\R^2$, since the foliations are 
non-singular (Poincare-Bendixson). Furthermore the leafs of $\widetilde{\mathfrak{F}}^\pm$ through $x\in \widetilde{M}$ 
lies at a bounded distance from $x+m^\pm:=\{y\in\widetilde{M}|\;y-x\in m^\pm\}$ for any $x\in\widetilde{M}$. 

Consider a $k\in H_1(M^2,\Z)\cong \pi_1(M^2)$ with $k\in \inte\conv(\overline{m}^+\cup\overline{m}^-)$. 
Since $m^+\neq m^-$ the halflines $\overline{m}^+$ and $k-\overline{m}^-$ intersect exactly once. 
But then the future pointing lightlike geodesic through $x$ and tangential to $\widetilde{\mathfrak{F}}^+$ intersects the past 
pointing lightlike geodesic through $x+k$ and tangential to $\widetilde{\mathfrak{F}}^-$ exactly once. The resulting curve is 
future pointing and not a lightlike geodesic. Therefore $x+k\in I^+(x)$ for all $x\in \widetilde{M}$. This property implies
the viciousness of $(M,g)$. 

(ii)$\Rightarrow$(iii): If $(M,-g)$ is vicious there exists a pair of generators $\{\sigma,\tau\}$ of $\pi_1(M)$ 
which are represented by timelike loops. If the prime fundamental class $\sigma$ of every timelike loop is unique, up to
inversion, we can choose $\tau'\in\pi_1(M)$ such that $\sigma,\tau'$ generate $\pi_1(M)$. The viciousness of $(M,-g)$ 
implies that every point lies on a timelike loop. It is then easy to construct a timelike loop representing $\sigma^k\ast
\tau'$ for some $k\in\Z$. This contradicts the uniqueness of $\sigma$. 

To show that the universal cover of $(M,g)$ is globally hyperbolic is now easy. First every Lorentzian metric 
on a simply connected surface is causal. Second, since timelike and spacelike loops in $M$ with prime fundamental class 
can intersect at most once, 
every pair of lifts to $\widetilde{M}$ intersects at most once. Since the set of fundamental classes represented by 
spacelike loops generates $\pi_1(M)$, this shows that the sets $J^+(x)\cap J^-(q)$ are bounded (It is clear that they are 
closed). 

(iii)$\Rightarrow$(i): Since $(M,g)$ is vicious every point lies on a timelike loop. 
Consequently $\overline{m}^+$ and $\overline{m}^-$ cannot coincide. If we assume $m^+=m^-$ then $\overline{m}^+=-\overline{m}^-$. 
In this case the causal diamonds $J^+(x)\cap J^-(y)$ in the universal cover cannot be compact. 
In fact since future pointing lightlike curves tangential to different lightlike distributions in $M$ run in, up to a constant,
opposing directions, it is easy, using the viciousness of $(M,g)$, to find a constant $C<\infty$ such that for any 
pair of causally related points $y\in J^+(x)\subset\widetilde{M}$ of Riemannian distance at least $C$ and every 
$K>0$ there exits $z\in J^+(x)\cap J^-(y)$ and $\dist(x,z),\dist(y,z)\ge K$.
\end{proof}

We obtain the corollaries.

\begin{cor}\label{C60}
A closed Lorentzian surface $(M,g)$ is class $A$ if and only if $(M,-g)$ is class $A$.
\end{cor}

\begin{cor}\label{C60a}
The universal cover $(\widetilde{M},\widetilde{g})$ of a class A surface $(M,g)$ is geodesically connected.
\end{cor}

\begin{proof}[Proof of corollary \ref{C60a}]
Since $(\widetilde{M},\widetilde{g})$ as well as $(\widetilde{M},-\widetilde{g})$ are globally hyperbolic, it remains to note that every pair of points 
$(p,q)\in \widetilde{M}\times\widetilde{M}$ is causally related for $(\widetilde{M},\widetilde{g})$ or $(\widetilde{M},-\widetilde{g})$. But this fact follows 
from the observation that the leaf of $\widetilde{\mathfrak{F}}^+$ through $p$ intersects the leaf of $\widetilde{\mathfrak{F}}^-$ through $q$ exactly once, if 
$(M,g)$ is class A.
\end{proof}

\begin{remark}
For closed Lorentzian surfaces $(M,g)$ such that the fundamental classes of all causally constant loops belong to a subgroup $\langle \sigma \rangle <\pi_1(M)$
generated by $\sigma \in \pi_1(M)$, the universal cover of $(M,g)$ is not geodesically connected. This condition is satisfied for all Lorentzian surfaces with 
$k_g\neq 0$.
\end{remark}

\begin{prop}
All Lorentzian metrics $g$ with $k_g>0$ are of class $B$.
\end{prop}

\begin{proof}
If $m^+\neq m^-$, proposition \ref{P49} implies $k_g=0$.
\end{proof}

\subsection{Closed geodesics in class A surfaces}

In the case of class $A$ surfaces we can refine theorem \ref{T60} to the following result. 
\begin{theorem}\label{T60b}
Let $(M,g)$ be of class $A$. Then every free homotopy class contains a closed geodesic. 
More precisely, if the free homotopy class contains timelike (spacelike) loops only, then all closed geodesics 
in this class are timelike (spacelike). If a free homotopy class is represented by time- and spacelike loops, then the class 
contains at least one timelike and one spacelike closed geodesic. If the free homotopy class is represented by lightlike loops 
only, then all causally constant loops (up to parameterization) in that class are contained in that foliation and are closed lightlike geodesics.
\end{theorem}

Note that parts of theorem \ref{T60b} were already observed in lemma \ref{A2}.

Before we proceed to prove theorem \ref{T60b}, we want to make some comments on how this result generalizes known ones 
about conformally stationary Lorentzian surfaces. We have the following simple proposition.

\begin{prop}
Every closed conformally stationary Lorentzian surface is class A.
\end{prop}

\begin{proof}
We can assume w.l.o.g. that $(M,g)$ is time and space orientable. If $(M,g)$ is class B then it contains a compact leaf of either lightlike foliation. The conformal 
Killing vector field $X$ has to be transversal to this compact leaf since $X$ is assumed to be timelike. Further the flow of $X$ maps compact leafs to compact leafs. 
If the rotation class of $X$ differs from $m^\pm$, then every flowline of $X$ intersects every leaf of $\mathfrak{F}^+$ and $\mathfrak{F}^-$. But in this case both 
lightlike foliations have to consist entirely of compact leafs. We saw earlier that in this case $(M,g)$ is class A.

If the rotation class of $X$ is equal to $m^\pm$ and $(M,g)$ is class B then the flow of $X$ has a periodic orbit and all orbits converge to a periodic one. Under 
these circumstances and together with the observation that compact leafs of $\mathfrak{F}^\pm$ are mapped to compact ones by the flow of $X$, we see that there 
exists a compact leaf of $\mathfrak{F}^\pm$ which intersects a periodic orbit of $X$ twice and has the same fundamental class. This is impossible. Therefore 
$(M,g)$ has to be class A.
\end{proof}

Theorem \ref{T60b} now implies.
\begin{cor}
Every conformally stationary closed Lorentzian surface contains infinitely many closed time and spacelike geodesics.
\end{cor}

This proof of the existence of a closed geodesic in a conformally stationary closed surfaces is limited to the case $\dim M=2$. The case $\dim M\ge 3$ will 
require completely different methods. For $\dim M\ge 3$ the compactness argument is no longer valid. 

Now we turn to the proof of theorem \ref{T60b}.

\begin{lemma}
Every fundamental class of a class A surface is representable by a causally constant loop.
\end{lemma}

\begin{proof}
Choose time and space orientations for the universal cover. Then we orient the lightlike foliations $\widetilde{\mathfrak{F}}^\pm$ to be future pointing. 
Note that any pair of leafs not belonging to the same lightlike foliation intersect 
exactly once. This is due to the fact that the leafs of $\widetilde{\mathfrak{F}}^\pm$ through $x\in\widetilde{M}$ 
lie in a bounded neighborhood of the lines $x+m^+$ and $x+m^-$.

Let $\pi$ be a fundamental class of $M$. 
The leaf of $\widetilde{\mathfrak{F}}^+$ through $x$ intersects the leave of $\widetilde{\mathfrak{F}}^-$ through 
$x+\pi$ in $y$. If $y\in J^+(x)\cap J^-(x+\pi)$ or $y\in 
J^-(x)\cap J^+(x+\pi)$ the resulting curve $\gamma_{x}$ connecting $x$ with $x+\pi$ 
will be nonspacelike. If $y\in J^+(x+\pi)\cap J^+(x)$ or $y\in J^-(x+\pi)\cap J^-(x)$ the resulting curve 
$\gamma_{x}$ will be nontimelike.
\end{proof}

\begin{lemma}\label{L60}
Let $(M,g)$ be a Lorentzian surface and $Z\subset M$ an annulus bounded by lightlike curves and such that $(Z,g|_Z)$ contains a timelike or spacelike loop. Then 
$Z$ contains a closed maximal timelike or spacelike geodesic. If the interior of $Z$ does not contain any closed smooth lightlike curves the ``or'' is exclusive, i.e. all 
closed geodesics of $(\inte Z,g|_{\inte Z})$ are either timelike or spacelike.
\end{lemma}
Note that theorem \ref{schell10} is a special case of the claim. Observe that no assumption is made  in lemma \ref{L60} whether the interior of $Z$ contains closed 
leafs of $\mathfrak{F}^\pm$.

\begin{proof}
The first part will follow from the second assertion since any annulus bounded by lightlike
curves that contains at least one causally constant and nonlightlike loop has to contain an annulus bounded by 
lightlike curves and the interior does not contain any other closed smooth lightlike curves. 
Recall that under this assumption the only causally constant loops in $(Z,g|_Z)$ are either all nonspacelike or 
nontimelike.

Now consider the space $\Lambda_1(Z,[g|_Z])$ of loops in $Z$ of constant causal character with prim fundamental 
class in $\pi_1(Z)$. All loops in $\Lambda_1(Z,[g|_Z])$ can be given a Lipschitz parameterization. 
Therefore the Riemannian length functional is well defined on $\Lambda_1(Z,[g|_Z])$.
Note that every curve in $\Lambda_1(Z,[g|_Z])$ is free of self-intersections. Indeed if $\gamma\in \Lambda_1(Z,[g|_Z])$ 
has self-intersections, 
consider a lift $\widetilde{\gamma}$ of $\gamma$ to the universal cover of $Z$ and a translate $\widetilde{\gamma} 
+kv$ intersecting $\widetilde{\gamma}$ (wlog $k\ge 0$). Now construct a closed curve $\widetilde{\beta}$ in the 
universal cover as follows: Follow $\widetilde{\gamma}+kv$ from the starting point to an intersection with 
$\widetilde{\gamma}$. Change to $\widetilde{\gamma}$ and follow $\widetilde{\gamma}$ until the endpoint. To get back to 
the starting point just follow $\widetilde{\gamma}+k'v$ for each $1\le k'\le k$. This curve has, by assumption, 
constant causal character and is closed in the universal cover. Therefore the universal cover of $Z$ contains a 
nullhomotopic loop with constant causal character. This is impossible by the theorem of Poincare-Bendixson 
(\cite{hirschsmale}).  

To prove that $\Lambda_1(Z,[g|_Z])$ contains a closed geodesic, we show that the subset of $\Lambda_1(Z,[g|_Z])$ consisting 
of curves parameterized w.r.t. constant $g_R$-arclength is compact relative to the topology of uniform convergence, 
i.e. the Riemannian length functional is bounded on that subset.

Assume $\sup_{\gamma\in \Lambda_1(Z,[g|_Z])}\{L^{g_R}(\gamma)\}=\infty$. Extend $Z$ by two collars to a larger annulus $Z'$. 
Next choose a $g'\in Lor(Z)$ strictly wider than $g$ and extend $g'$ to $Z'$ such that $\partial Z'$ consists of 
compact leafs of the lightlike foliations of $(Z',g')$. Consider a sequence $\{\gamma_n\}\in \Lambda_1(Z,[g|_Z])$ 
with $L^{g_R}(\gamma_n)\to\infty$. W.l.o.g. we can assume that the $\gamma_n$ are smooth. Reparameterize each $\gamma_n$ 
w.r.t. $g_R$-arclength and consider the Borel propability measures 
$$\mu_n:=\frac{1}{L^{g_R}(\gamma_n)}(\frac{d}{ds}\gamma_n)_\sharp (\mathcal{L}^1)$$
on $TZ$, where $s$ is the $g_R$-arclength parameter. The sequence $\{\mu_n\}_{n\in\N}$ is bounded in $C^0(TZ,\R)'$ 
and therefore contains a converging subsequence. Choose a point $w\in TZ$ in the support of the limit measure. 
$w$ is a $g$-causal vector and consequently $g'$-timelike. But then we can construct a closed timelike curve in 
$(Z',g')$ with selfintersection. First note that $\pi_{TZ}(w)\in Z\subset \inte Z'$. Choose a convex neighborhood 
$U$ of $\pi_{TZ}(w)$ in $(Z',g')$ and a sequence $s_k\in \R$ such that $\gamma_{n_k}(s_k)\to \pi_{TZ}(w)$ for some 
sequence $n_k\in\N$. Since
$$\mathcal{L}^1((\frac{d}{ds}\gamma_n)^{-1}(TU))\ge \e L^{g_R}(\gamma_n)$$
for some $\e=\e(U)>0$, $\gamma_n$ has to leave $U$ and return more than once arbitrarily close to $\pi_{TZ}(w)$. 
Therefore we can choose parameter values 
$a_k,b_k$ belonging to different connected components of $\gamma_n^{-1}(U)$ such that $\gamma_{n_k}(a_k),
\gamma_{n_k}(b_k)\to \pi_{TZ}(w)$. Since all tangents of $\gamma_n$ are equi-timelike relative to $g'$ in the sense 
that $\dist(\frac{d}{ds}\gamma_n,\Light(Z',[g']))\ge\e$ for some $\e>0$, we can deform $\gamma_n$ in $U$ to intersect 
itself for some $n$ sufficiently large. This is impossible by the first part of the proof. 
Therefore we have $\sup_{\gamma\in \Lambda_1(Z,[g|_Z])}\{L^{g_R}(\gamma)\}<\infty$. 

Theorem \ref{T39} then 
ensures the existence of a closed maximal geodesic in $(Z,g|_Z)$. Note again that with $(M^2,g)$ Lorentzian, 
$(M^2,-g)$ is Lorentzian as well, so theorem \ref{T39} is applicable.
\end{proof}

\begin{proof}[Proof of theorem \ref{T60b}]
We can divide the homotopy classes $\pi$ representable by causally constant curves into four sets. Those for which 
$\Lambda(M,[g])_\pi$ consits solely of timelike loops, those such that $\Lambda(M,[g])_\pi$ consists of spacelike loops, 
those such that $\Lambda(M,[g])_\pi$ consists of lightlike loops and those such that $\Lambda(M,[g])_\pi$
contains timelike as well as spacelike loops. 

By corollary \ref{C60} the second case follows from the first. The first case has been treated essentially
in \cite{schelling}. Lemma \ref{A2} and corollary \ref{CD10} in appendix \ref{AB} show that every integer homology class in the 
time-orientation cover of a class A surface that is contained in the interior of the convex hull of $\overline{m}^+\cup 
\overline{m}^-$ contains a closed future pointing timelike geodesic. Since the assumption that a homotopy class that gives 
rise to a foliation of the time orientation cover by future pointing timelike curves is equivalent to the image in $H_1(M,\R)$ 
being contained in that interior, the first case follows.

The case that the free homotopy class contains timelike as well as spacelike loops follows from lemma \ref{L60}. 
Because in that case the free homotopy class contains at least two disjoint closed lighlike loops. Any pair of disjoint 
lightlike loops then bounds an annulus in $M$. Thus the assumptions of lemma \ref{L60} are satisfied.

At last we have to prove that if the free homotopy class is represented by lightlike loops only, then all these loops are closed lightlike geodesics. First it is easy to 
see that all lightlike geodesics parameterizing the lightlike loops are periodic. Next it follows that these periodic geodesics form a foliation of $M$, since $g$ is 
smooth. 

Now assume that one periodic geodesic $\gamma$ in that foliation is incomplete. Then \cite{hawk} proposition 6.4.4 implies that there exists a timelike loop freely 
homotopic to $\gamma$. This contradicts our assumption that all causally constant loops in that free homotopy class are lightlike.  Therefore all the periodic 
geodesics in the foliation are closed.
\end{proof}

\subsection{Closed geodesics in class B surfaces}\label{lastB}

In the case of $m^+=m^-$ no leaf of either one of the 
lightlike foliations can intersect all leafs of the other foliation. Therefore both foliations contain a compact leaf and 
$m^+,\; m^-$ are rational. Additionally, since the lightlike distributions are transversal, neither lightlike foliation can 
consist entirely of closed orbits. Furthermore any pair of smooth closed lightlike curves must be disjoint. Therefore we 
conclude that there exists a pair of closed leafs $\zeta_{1,2}$ of the lightlike foliations in $(M,g)$ bounding an 
annulus $Z$ in $M$ such that $(\inte Z,g_{\inte Z})$ contains no closed lightlike curves. Note that $(Z,g_Z)$ is time and 
space orientable. This follwos directly from that fact that $Z$ is orientable and the boundary curves are lightlike, i.e. 
causally constant.
Therefore lemma \ref{L60} implies theorem \ref{T60} for class B surfaces. Together with theorem \ref{T60b} this implies theorem \ref{T60} 
for all closed Lorentzian surfaces.

\begin{proof}[Proof of theorem \ref{T60a}]
Theorem \ref{T60} ensures the existence of at least one closed geodesic in every compact Lorentzian surface.
If $(M,g)$ contains only finitely many closed geodesics then, by theorem \ref{T60b},  $m^+$ and $m^-$ must coincide. 

Consider the number of 
annuli bounded by closed smooth lightlike curves in $(M,g)$. If this number is at least two, then lemma \ref{L60} implies 
the existence of at least two closed geodesics. Therefore the only interesting case is if the number of annuli is 
equal to one. In this case $(M,g)$ contains exactly one closed lightlike curve $\gamma$ (up to parameterization). 
The proof consists now of showing that $\gamma$ can be parameterized as a closed lightlike geodesic.

By switching from $g$ to $-g$ we can assume that $(M,g)$ contains only non-spacelike closed curves. Then there 
exists an $\e>0$, a neighborhood $U$ of $\gamma$ and a smooth function $f\colon S^1\times (-\e,\e)\to \R$ with 
$f|_{S^1\times\{0\}}\equiv 0$ such that $(U,g|_U)$ is isometric to the Lorentzian annulus 
$(S^1\times (-\e,\e),g')$, where $g'$ is expressed in the natural coordinates $(\phi,s)$ as 
$g'=-d\phi (ds-f(\phi,s)d\phi)$. The lightlike distributions of $(S^1\times (-\e,\e),g')$ are spanned by 
$\partial_s$ and $\partial_\phi +f(\phi,s)\partial_s$. Then the property that $(S^1\times (-\e,\e),g')$ contains only 
nonspacelike closed curves can be expressed as $\int_0^T f(\phi(\tau),s(\tau))d\tau \le 0$, where 
$\tau\mapsto (\phi(\tau),s(\tau))$ is any curve tangential to $\partial_\phi +f(\phi,s)\partial_s$ with 
$\phi(T)=\phi(0)$ and $\dot{\phi}\ge 0$.

The $\phi$-coordinate part of the equation of the geodesic is $\ddot{\phi}+\partial_s f \dot{\phi}^2=0$.
Consider a solution $\tau\mapsto (\phi(\tau),0)$ with $\phi(T)=\phi(0)$.
Then the lightlike geodesic parameterizing $\{s=0\}$ is closed if and only if $\int_0^T \partial_s f(\phi(\tau),0)d\tau=0$. 
This integral equation is an easy consequence of the above integral inequality.
\end{proof}

A closer look reveals that the least number of closed geodesics in 
$(M,g)$ can be connected to the homotopy class of $g$ in $Lor(M)$. 
Recall the definition of $k_g:=\min\{|n_\sigma(g)|+|n_\tau(g)|\}$, where $g$ is time and space orientable Lorentzian metric 
on $M$ and $\sigma,\tau$ generate $\pi_1(M)$.
\begin{prop}\label{T61}
Let $(M,g)$ be a space and time orientable closed Lorentzian surface. Then $(M,g)$ contains
\begin{enumerate}
\item  $4\cdot k_g$ closed leafs of the lightlike foliations and
\item $2\cdot k_g$ closed timelike as well as $2\cdot k_g$ closed spacelike geodesics.
\end{enumerate}
\end{prop}

For $(M,g)$ not time and space orientable we set $k_g:=k_{g'}$, where $(M',g')$ is a time and 
space orientable minimal cover. Minimal is the sense that there is no time and space orientable cover
with smaller deck transformation group. Then we obtain the following corollary, via simple counting argument.

\begin{cor}
Every closed Lorentzian surface $(M,g)$ contains at least $k_g$-many timelike as well as spacelike 
closed geodesics. 
\end{cor}

\begin{exmp}\label{E60}
Explicit examples show that this lower bound on the multiplicity is optimal. Consider for $k\in \Z$ on $\mathbb R^2$ 
the Lorentzian metrics: 
    $$\overline{g}_k(x,y):= \sin(4\pi k x)[dx^2-dy^2]+\cos(4\pi k x)[dxdy]$$
These metrics induce Lorentzian metrics $g_k$ on $T^2=\R^2/\Z^2$ with $k_{g_k}=|k|$.
Then the equation of the geodesic for $\overline{g}_k$ reads
\begin{align*}
\ddot{x}&-2\pi k \sin(4\pi kx)\cos(4\pi kx) \dot{x}^2\\
   &+4\pi k\cos^2(4\pi kx) \dot{x}\dot{y}+2\pi k  \cos(4\pi kx)\sin(4\pi kx) \dot{y}^2=0\\
\ddot{y}&-2 \pi k(1+\sin^2(4\pi kx))\dot{x}^2\\
   &+4\pi k\sin(4\pi kx) \cos(4\pi kx) \dot{x}\dot{y}-2 \pi k\cos^2(4\pi kx) \dot{y}^2=0.
\end{align*}
Recall that closed causally constant curves cannot leave strips bounded by lightlike curves. 
Consequently any closed geodesic must be contained in one of the strips $x^{-1}(\frac{l}{k}+\frac{i}{4k},\frac{l}{k}
+\frac{i+1}{4k}),\; 0\le l\le k-1,\;0\le i\le 3$. 
If $\gamma\colon (\alpha_\gamma,\omega_\gamma)\to \R^2$, $t\mapsto \gamma(t)=:(x(t),y(t))$ is a geodesic with 
$\dot{x}(t_0)=0$ $(t_0\in (\alpha_\gamma,\omega_\gamma))$ in $x^{-1}(\frac{l}{k},\frac{l}{k}+\frac{1}{4k})$ or 
$x^{-1}(\frac{l}{k}+\frac{1}{2k},\frac{l}{k}+\frac{3}{4k})$ $(0\le l\le k-1)$, then
$$\ddot{x}(t_0)=-2\pi k\cos(4\pi kx(t_0)) \sin(4\pi kx(t_0)) \dot{y}^2(t_0)<0.$$ 
Therefore the closed geodesics in these strips have to be constant in the $x$-coordinate. 
Note that in the present sign convention these closed geodesics are timelike. The timelike closed geodesics parameterize 
the sets $x^{-1}(\frac{8l+1}{8k})$ and $x^{-1}(\frac{8l+5}{8k})$.
The same argument, only with signs reversed, applies to the closed (spacelike) geodesics in 
$x^{-1}(\frac{l}{k}+\frac{1}{4k},\frac{l}{k}+\frac{1}{2k})$ and $x^{-1}(\frac{l}{k}+\frac{3}{4k},\frac{l+1}{k})$. This 
shows that $(T^2,g_k)$ contains exactly $2\cdot k_g$ many closed timelike and $2\cdot k_g$ many closed spacelike 
geodesics. The number of closed leafs of the lightlike foliations is exactly $4\cdot k_g$.
The geodesic flow of these manifolds is discussed in greater generality in \cite{sa1}, Appendix, case (2C)(b)).
\end{exmp}

\begin{proof}[Proof of proposition \ref{T61}]
The claim is trivial for $k_g=0$. Therefore we can assume that $k_g>0$.

(1) Consider the two transversal lightlike oriented foliations $\mathfrak{F}^+,\mathfrak{F}^-$ of $(M,g)$. 
We know that $\mathfrak{F}^+$ or $\mathfrak{F}^-$ contain closed leafs, since $k_g>0$.

Choose a prime fundamental class $\pi\in\pi_1(M)$ whose image spans the subspace representing $m^+=m^-$. 
Then all simply closed lightlike curves in $(M,g)$ have fundamental class either $\pi$ or $\pi^{-1}$. 
Divide the periodic orbits of $\mathfrak{F}^+$ into $\{\zeta^+_\alpha\}_{\alpha\in A}$ such that 
$[\zeta^+_\alpha]=\pi$ and $\{\eta^-_\beta\}_{\beta\in B}$ s.t. $[\eta^-_\beta]=\pi^{-1}$. The assertion will follow 
from: $\sharp A$, $\sharp B\ge k_g$. 

Choose a vector field $X^+\in \Gamma(TM)$ tangential to $\mathfrak{F}^+$ and compatible with the orientation of 
$\mathfrak{F}^+$. Assume to the contrary that $\sharp A <k_g$. Cut $M$ along the $\zeta^+_\alpha$ in $\sharp A$-many 
annuli. These annuli contain, except on the boundary, no closed orbits of $\mathfrak{F}^+$ with fundamental class $\pi$. 
Let $Z$ be one of the annuli defined by $\mathfrak{F}^+$. Then $X^+$ restricts to a vector field on $Z$ tangential to 
the boundary. Identifying the boundary curves pointwise gives a torus $T^2$. $X^+|_Z$ projects to a vector field $X'$ on 
$T^2$. Choose a fundamental class $[\alpha]\in \pi_1(T^2)$ not contained in the image of $\pi_1(Z)$ under the 
projection. Now define the rotation number of $\mathfrak{F}^+|_Z$ as the degree of $X'\circ \alpha$. Note that any curve 
in $Z$ connecting the boundary components can be extended along the boundary such that the projection to $T^2$ of the 
extension is closed. The fundamental class of the projection is not in the image of $\pi_1(Z)$. Therefore the 
rotation number can be defined for curves connecting the boundary components of $Z$. The claim is that the absolute 
value of the rotation number on any annulus is at most one. Assume the rotation number to be non-vanishing. Then the 
number of Reeb-components in $Z$ must be even and at most two, since the boundaries of Reeb-components represent 
inverse homotopy classes in $\pi_1(Z)$. Note that the boundary must be part of the limit of these Reeb-components.
Denote with $Z'$ the annulus bounded by the limit curves of the Reeb-components that are not the boundary. Since
$\mathfrak{F}^+|_{Z'}$ contains no Reeb-components, $Z'$ gives rise to a regular curve $\gamma_{Z'}$ connecting the
boundary components of $Z'$ and transversal to $\mathfrak{F}^+|_{Z'}$. Extend $\gamma_{Z'}$ to a curve $\gamma_Z$ 
connecting the boundary components of $Z$ such that the additional arcs do not intersect $Z'$.  Since $TZ$ is trivial 
one can choose a vector field $Y_Z$ along $\gamma_Z$ equal to $X^+$ on $\gamma_{Z}|_{\partial Z\cup Z'}$ and transversal 
to $\dot{\gamma}_Z$ on $\gamma_Z\setminus \gamma_{Z'}$. The rotation number of $Y_Z$ along $\gamma_Z$ 
obviously vanishes. The rotation of $X^+$ along $\gamma_Z$ in one of the Reeb-components is $\pm 1/2$. Consequently 
the entire rotation number can only be $\pm 1$ (compare \cite{hehi} theorem 4.2.15). 

Choose any fundamental class $\pi'$ such that $\{\pi,\pi'\}$ generate $H_1(M,\Z)_\R$. 
Let $\gamma\colon S^1\to M$ be any representative of $\pi'$. By deforming $\gamma$ suitably we can 
assume that the restriction of $\gamma$ to every annulus connects the boundary components of this annulus. It is easy 
to see that the sum of the rotation numbers of the restrictions of $\gamma$ to the annuli is $k_g$.  Since the rotation 
number on each annulus is $\pm 1$ or $0$, the sum is strictly smaller than $k_g$. This contradicts the assumption. 
Therefore $\sharp A, \sharp B\ge k_g$. The claim now follows since $(M,g)$ contains two transversal lightlike 
foliations.

(2) Consider the annuli in $M$ bounded by closed leafs of $\mathfrak{F}^+$ and containing exactly one 
Reeb-component of $\mathfrak{F}^+$.
Let $Z$ be one of these annuli. Choose vector fields $X^+$ tangential to $\mathfrak{F}^+$ and $X^-$ tangential to
$\mathfrak{F}^-$ such that $X^-|_{\partial Z}$ points ``outwards'' and $\conv\{X^+,X^-\}$ consists of timelike vectors. 
Note that the orbits of $X^+$ parameterizing $\partial Z$ represent inverse fundamental classes. 
Consider the $\alpha$-limit of the $X^-$-orbits through points in $\partial Z$. These $\alpha$-limits are compact leafs 
of $\mathfrak{F}^-$ and must be contained in $Z$ (Reeb-components do not have transversals).
The boundary curves of $Z$ together with the associated $\alpha$-limits form two sub-annuli  $Z^\pm$. Labels are 
distributed as follows: Denote $Z^-$ the annulus whose boundary orbits, one of $X^+$ and one of $X^-$, share the same 
fundamental class. Exactly one of the annuli satisfies this condition. In fact the $\alpha$-limit of the orbits of $X^+$ 
through points in $Z^\circ$ is the boundary curve of $Z^-$ belonging to $\mathfrak{F}^+$. The other annulus 
will be denoted by $Z^+$. Note that on $Z^-$ the lightlike vector fields which are not tangential to the boundary point 
outwards. Therefore $Z^-$ contains points whose $\omega$-limit of both $X^+$ and $X^-$ do not belong to $Z^-$. 
This shows that $Z^-$ contains a closed timelike curve. Lemma \ref{L60} then implies the existence of a closed timelike 
geodesic in $Z^-$. The same argument applies to $Z^+$ except that $X^+$ has to be replaced by $-X^+$. This shows 
that $Z^+$ contains a closed spacelike geodesic. Part (1) showed that every lightlike foliation contains at least 
$2\cdot k_g$-many Reeb components. This implies part (2) of the proposition.
\end{proof}

\appendix

\section{Schelling's Diploma Thesis}\label{AB}

The results of \cite{schelling} represent important special cases of the results in the present text. For that reason and since they are not published anywhere, we 
will recollect the relevant parts in this appendix.

The diploma thesis \cite{schelling} considers closed Lorentzian $2$-manifolds $(M,g)$ such that both lightlike distributions are orientable. We have seen earlier 
that this implies $M$ to be orientable and therefore $M$ has to be diffeomorphic to the $2$-torus. First we review the results on the case $m^+\neq m^-$. 

\begin{theorem}[\cite{schelling}, Satz 4.1]\label{A1}
Let $(M,g)$ be a $2$-dimensional closed oriented spacetime such that $m^+\neq m^-$. Then the Abelian cover is globally hyperbolic.
\end{theorem}
The claim is contained in proposition \ref{P60}. Therefore we skip the proof.

\begin{lemma}[\cite{schelling}, Lemma 4.3]\label{A2}
Let $h\in \inte \conv(\overline{m}^+\cup \overline{m}^-)\cap H_1(M,\Z)$. Then there exists a closed timelike geodesic $\gamma\colon \R\to M$ representing 
$h$ and maximizing arclength among all causal curves with homology class $h$.
\end{lemma}

The proof is an adaptation of well known techniques from \cite{hed1}. It includes a maximization argument on the space of closed causal curves re\-pre\-sen\-ting 
$h$ very much like the argument for Tipler's theorem. In fact it is possible to reduce the claim to Tipler's theorem: Choose any simply connected smooth spacelike 
curve $\eta$ in $(M,g)$. Since $(M,g)$ is vicious (proposition \ref{P60}) $\eta$ is not acausal in $(M,g)$. Further the proof of proposition \ref{P60} directly implies 
that the Geroch covering of $M$ relative to $\eta$ is globally hyperbolic. Now Tipler's theorem implies the claim. It should be noted that this is not the argument 
employed in \cite{schelling}. Therein the claim is proved directly, repeating Tipler's argument in this special case.

As in \cite{hed1}, the fact that we are considering geodesics in $2$-manifolds, gives further information on the minimal period of the closed geodesic:

\begin{cor}[\cite{schelling}, Korollar 4.4]\label{CD10}
Let $\gamma\colon \R\to M$ be a closed timelike geodesic with homology class $h$ maximizing arclength among all causal representatives of $h$. 
Denote with $T$ the minimal period of $\gamma$. Then the class $h$ is relative prim in $H_1(M,\Z)$, i.e. for any $h'\in H_1(M,\Z)$ and 
$\lambda >0$ with $h=\lambda h'$ we have $\lambda =1$ and $h'=h$.
\end{cor}

The proof relies on the observation that if $\gamma\ast\gamma$ is not maximal among all causal curves representing $2h$, then the closed timelike
geodesic maximizing the arclength for $2h$ will intersect $\gamma\ast\gamma$ at least twice. These intersection are transversal since $\gamma\ast\gamma$ is 
not maximal. From this we could construct a causal curve representing $h$ with greater arclength that $\gamma$. This is a contradiction to our assumption.

\cite{schelling} contains partial results for the case $m^+ = m^-$ as well. In \cite{schelling} the approach is different from the one in this article, insofar as 
\cite{schelling} focuses on the dynamics of the foliations $\mathfrak{F}^\pm$. Our approach, on the other hand, relies rather on the set of causally constant loops. 

Two disjoint closed leafs $\zeta_{1,2}$ of $\mathfrak{F}^\pm$ are called {\it neighboring} if the interior of the annulus bounded by $\zeta_{1,2}$ contains no closed 
leaf of $\mathfrak{F}^+$ or $\mathfrak{F}^-$. This notion naturally extends to the universal covering, by referring to periodic leafs rather than closed ones.

\begin{lemma}[\cite{schelling} Lemma 5.2]
Let $\zeta^\pm$ be neighboring leaves of $\mathfrak{F}^\pm$. Denote $B$ the strip bounded by $\zeta^\pm$. Let $p\in B^\circ\cup \zeta^-$ and $q\in B^\circ 
\cup \zeta^+$. Denote $\zeta^+_p$ the leaf of $\mathfrak{F}^+$ through $p$ and $\zeta^-_q$ the leaf of $\mathfrak{F}^-$ through $q$. Then $\zeta^+_p$ and 
$\zeta^-_q$ intersect exactly once. Further the intersection lies in $B^\circ$.
\end{lemma}

As we have seen above, the matter of finding maximal timelike geodesics in the case $m^+ = m^-$ is more subtle than in the case of $m^+ \neq m^-$. In 
\cite{schelling} a criterion is given whenever between two neighboring leafs of $\mathfrak{D}^\pm$ a maximal timelike periodic geodesics exists. The results of 
section \ref{lastB} extend these special cases to the general case.

\begin{definition}[\cite{schelling}]
Let $m^+=m^-$. A non-periodic future pointing lightlike geodesic $\eta$ is {\it affiliated} to a periodic leaf $\zeta$ of $\mathfrak{F}^\pm$, if 
$\eta$ is either $\omega$-asymptotic towards $\zeta$ or $\eta$ intersects $\zeta$.
\end{definition}

Denote by $\eta^\pm_p$ the future pointing geodesic parameterization of the leaf of $\mathfrak{F}^\pm$ through $p\in M$ with $\eta^\pm_p(0)=p$.

\begin{theorem}[\cite{schelling}, Satz 6.1]\label{schell10}
Let $m^+=m^-$ and $\zeta_1$ and $\zeta_2$ be neighboring leafs. Denote by $B$ the connected component between $\zeta_1$ and $\zeta_2$. Let $p\in B^\circ$. 
Then the following assertions are equivalent:
\begin{enumerate}
\item $\eta^+_p$ and $\eta^-_p$ are affiliated to different $\zeta_1$ and $\zeta_2$.
\item $B$ contains a maximal periodic timelike geodesic.
\end{enumerate}
\end{theorem}

\bibliographystyle{amsalpha}

\end{document}